# GLOBAL STABILIZATION BY MEANS OF DISCRETE-DELAY OUTPUT FEEDBACK

**Iasson Karafyllis**

**Department of Environmental Engineering, Technical University of Crete,
73100, Chania, Greece
email: ikarafyl@enveng.tuc.gr**

**Abstract**
In the present work, sufficient conditions for global stabilization of nonlinear uncertain systems by means of discrete-delay static output feedback are presented. Illustrating examples show the efficiency of the proposed control strategy.

## 1. Introduction

Feedback laws with delayed terms have been considered recently for the solution of various control problems. Particularly, the following works have showed that feedback laws, which involve delays, present features that cannot be induced by means of ordinary feedback (i.e., feedback with no delays):

- in [1], it is shown that the use of discrete-delays can allow the design of observers that provide state estimates for linear systems which converge to the actual state values in a pre-specified finite settling time,
- in [6], it is shown that the use of distributed delay feedback can allow the design of smooth feedback laws which achieve finite-time stabilization of nonlinear systems in a pre-specified finite settling time,
- in [3,8,9,10,11], it is shown that the use of discrete delay static output feedback can achieve stabilization for linear systems, which cannot be stabilized by ordinary (i.e., with no delays) static output feedback

Specifically, in [3] it was shown that autonomous, minimum phase, linear systems $\dot{x} = Ax + Bu$, $x \in \Re^n$, $u \in \Re^m$, $y(t) = cx(t)$ with relative degree 1 or 2 can be stabilized by static output feedback with delays of the form $u(t) = -k_1 y(t) - k_2 \frac{y(t) - y(t-h)}{h}$, where $h > 0$. The form of the feedback is obtained by replacing the derivative $\dot{y}(t)$ in an ordinary stabilizing feedback $u(t) = -k_1 y(t) - k_2 \dot{y}(t)$ by the numerical approximation of the derivative of the output signal $\frac{y(t) - y(t-h)}{h} \approx \dot{y}(t)$. The same idea was used in [11] for the stabilization of a chain of integrators using multiple delays.

In the present work, we generalize the idea proposed in [3] and we show that for uncertain systems of the form:

$$\dot{x}_i(t) = x_{i+1}(t) + v_i(t) \quad , \quad i = 1,\ldots,n-1$$
$$\dot{x}_n = v_n(t) + a(d(t), x(t))u(t) \tag{1.1}$$
$$x(t) = (x_1(t),\ldots,x_n(t))' \in \Re^n, d(t) \in D, v(t) = (v_1(t),\ldots,v_n(t))' \in \Re^n, u(t) \in \Re$$

where $D \subseteq \Re^l$ is compact, $d = d(\cdot), v = v(\cdot)$ are time-varying uncertainties/disturbances (modelling errors), and $a : D \times \Re^n \to \Re$ is a locally Lipschitz function that satisfies the following inequality for certain constants $\alpha, \beta > 0$:

$$\alpha \leq a(d,x) \leq \beta \, , \text{ for all } (d,x) \in D \times \Re^n \tag{1.2}$$

there exists a vector $k \in \Re^n$ such that for sufficiently small $h > 0$, the linear feedback law with discrete delays:

$$u(t) = k' \begin{bmatrix} x_1(t) + e(t) \\ x_1(t-h) + e(t-h) \\ \vdots \\ x_1(t-(n-1)h) + e(t-(n-1)h) \end{bmatrix} \quad (1.3)$$

where $e : \Re \to \Re$ represents the measurement error, achieves robust global stabilization of the equilibrium point $0 \in \Re^n$ of (1.1) in the sense that the solution of (1.1) with (1.3) and arbitrary continuous initial condition $x_0 : [-(n-1)h, 0] \to \Re^n$ corresponding to arbitrary measurable and locally essentially bounded inputs $v : \Re^+ \to \Re^n$, $d : \Re^+ \to D$, $e : \Re \to \Re$, satisfies the following estimate for all $t \geq 0$:

$$|x(t)| \leq Q_0 \exp(-\mu t) \sup_{-(n-1)h \leq \theta \leq 0} |x(\theta)| + \sum_{i=1}^{n} Q_i \sup_{0 \leq \tau \leq t} \exp(-\mu(t-\tau))|v_i(\tau)| + Q_e \sup_{-(n-1)h \leq \tau \leq t} \exp(-\mu(t-\tau))|e(\tau)| \quad (1.4)$$

for appropriate constants $Q_0, Q_1, \ldots, Q_n, Q_e > 0$ (Proposition 2.3). Inequality (1.4) is a "fading memory estimate" that guarantees the Input-to-State Stability (ISS) property for the closed-loop system (1.1) with (1.3) (see [13] for the definition of the ISS property for finite-dimensional systems, [4,12] for the extension of the ISS property to time-delay systems and [5] for a discussion of "fading memory estimates" to a wide class of systems). The feedback law (1.3) is constructed using a backward difference numerical differentiator of the output signal $y(t) = x_1(t) + e(t)$, which is corrupted by measurement noise. Moreover, we obtain explicit estimates of the maximum allowable time step $h > 0$ that can be used for robust stabilization as well as of the gains $Q_1, \ldots, Q_n, Q_e > 0$ involved in (1.4). The obtained results are applied to minimum phase nonlinear systems of the form:

$$\begin{aligned} \dot{z}(t) &= f(d(t), z(t), x(t)), \; z(t) \in \Re^k \\ \dot{x}_i(t) &= g_i(d(t), z(t), x(t)) + x_{i+1}(t) + v_i(t) \quad, \quad i = 1, \ldots, n-1 \\ \dot{x}_n(t) &= g_n(d(t), z(t), x(t)) + a(d(t), z(t), x(t))u(t) + v_n(t) \\ y(t) &= x_1(t) + e(t) \\ x(t) &= (x_1(t), \ldots, x_n(t))' \in \Re^n, d(t) \in D, u(t) \in \Re \\ v(t) &= (v_1(t), \ldots, v_n(t))' \in \Re^n, e(t) \in \Re \end{aligned} \quad (1.5)$$

and we present sufficient conditions for robust global asymptotic stabilization by means of the discrete delay output feedback (1.3) (see hypotheses (A1-3) in the statement of Theorem 2.6). Robustness of the closed-loop system (1.5) with (1.3) is guaranteed by showing an inequality similar to (1.4) (Theorem 2.6). Therefore, the Bounded-Input-Bounded-State property and the Converging-Input-Converging-State property hold for the nonlinear closed-loop system (1.5) with (1.3). Moreover, the states $x(t)$ converge exponentially to the origin for the unforced case $v \equiv 0, e \equiv 0$. It should be noted that systems of the form (1.1) or (1.5) under the hypotheses imposed in the present work (see hypotheses (A1-3) in the statement of Theorem 2.6) can be stabilized by dynamic (observer-based) output feedback (see [14]). Consequently, it is clear that static output discrete-delay feedback is an alternative to dynamic output feedback and further studies need to be performed, in order to show the advantages and disadvantages of each type of feedback.

Consequently, the contribution of the paper is:

- the main result in [3] is generalized to uncertain minimum phase nonlinear systems of the form (1.5) with arbitrary relative degree (Theorem 2.6),
- the main result in [11] is generalized to uncertain systems of the form (1.1) and a completely different proof is provided (Proposition 2.3),
- robustness of the closed-loop system with respect to measurement and modelling errors is guaranteed by "fading memory estimates" of the form (1.4),
- explicit estimates of the maximum allowable time step $h > 0$ that can be used in (1.3) as well as of the gains of the inputs are provided.

The structure of the paper is as follows: in Section 2 the main results are presented and proved. Section 3 contains illustrating examples, which show the efficiency of the discrete-delay static output feedback for the stabilization of linear and nonlinear uncertain systems. Finally, in Section 4 we present the concluding remarks of the present work.

**Notations** Throughout this paper we adopt the following notations:

∗ By $Z^+$ we denote the set of non-negative integers and by $\Re^+$ we denote the set of non-negative real numbers.
∗ Let $I \subseteq \Re^+ := [0,+\infty)$ be an interval. By $L^\infty(I;U)$ ($L^\infty_{loc}(I;U)$) we denote the space of measurable and (locally) essentially bounded functions $u(\cdot)$ defined on $I$ and taking values in $U \subseteq \Re^m$. The sup operator used for a measurable and (locally) essentially bounded operator used in the text is actually the essential supremum of this function.
∗ By $C^0(A;\Omega)$, we denote the class of continuous functions on $A$, which take values in $\Omega$.
∗ Let $x:[a-r,b) \to \Re^n$ with $b > a > -\infty$ and $r > 0$. By $T_r(t)x$ we denote the "$r$-history" of $x$ at time $t \in [a,b)$, i.e., $T_r(t)x := x(t+\theta); \theta \in [-r,0]$.
∗ For a vector $x \in \Re^n$ we denote by $|x|$ its usual Euclidean norm, by $x'$ its transpose and by $|A| := \sup\{|Ax|; x \in \Re^n, |x|=1\}$ the induced norm of a matrix $A \in \Re^{m \times n}$. For $x \in C^0([-r,0];\Re^n)$ we define $\|x\|_r := \max_{\theta \in [-r,0]} |x(\theta)|$.
∗ For the definitions of the classes $K$ and $K_\infty$, see [7].
∗ By $A = diag(l_1,l_2,...,l_n)$ we mean that the matrix $A = \{a_{ij}; i=1,...,n, j=1,...,n\}$ is diagonal with $a_{ii} = l_i$, for $i = 1,...,n$.
∗ We say that a function $V \in C^0(\Re^n;\Re)$ is radially unbounded if the following property holds: "$V(x)$ is bounded if and only if $|x|$ is bounded".

## 2. Main Results

We start by presenting a preliminary result on the numerical differentiation of the output signal of an uncertain linear system. The following lemma shows that there is a family of backward difference operators (parameterized by the time step $h$) that provides state estimates for a linear uncertain observable system.

**Lemma 2.1** *Consider the following system:*

$$\begin{aligned} \dot{x}_i &= x_{i+1} + u_i, \ i=1,...,n-1 \\ \dot{x}_n &= u_n \\ x &= (x_1,...,x_n)' \in \Re^n, u = (u_1,...,u_n) \in \Re^n \end{aligned} \tag{2.1}$$

*with $n \geq 2$. There exist constants $K_0, K_1,..., K_n > 0$ and a family of linear continuous operators $\Delta^n_h : L^\infty_{loc}([-(n-1)h,0];\Re) \to \Re^n$ (parameterized by $h \in (0,1]$), defined by*

$$\Delta^n_h y := Q^n_h P^{-1} [y(0) \ \ y(-h) \ \ y(-2h) \ \ ... \ \ y(-(n-1)h)]', \ \forall y \in L^\infty_{loc}([-(n-1)h,0];\Re) \tag{2.2}$$

*where* $Q^n_h := diag\left(1, \dfrac{1!}{(-h)}, \dfrac{2!}{(-h)^2},..., \dfrac{(n-1)!}{(-h)^{n-1}}\right) \in \Re^{n \times n}$ *and* $P = \begin{bmatrix} 1 & 0 & 0 & ... & 0 \\ 1 & 1 & 1^2 & ... & 1^{n-1} \\ 1 & 2 & 2^2 & ... & 2^{n-1} \\ \vdots & \vdots & \vdots & ... & \vdots \\ 1 & n-1 & (n-1)^2 & ... & (n-1)^{n-1} \end{bmatrix} \in \Re^{n \times n}$, *such that for every $h \in (0,1]$, $u_i \in L^\infty_{loc}(\Re^+;\Re)$ ($i=1,...,n$), $x_0 \in \Re^n$ the solution $x(t)$ of (2.1) with initial condition $x(t_0) = x_0$ corresponding to inputs $u_i \in L^\infty_{loc}(\Re^+;\Re)$ ($i=1,...,n$), satisfies the following inequality:*

$$\left| x(t) - \Delta^n_h T_{(n-1)h}(t)x_1 \right| \leq \sum_{j=1}^n K_j h^{j+1-n} \sup_{t-(n-1)h \leq \tau \leq t} |u_j(\tau)|, \ \forall t \geq t_0 + (n-1)h \tag{2.3}$$

*Moreover, for every $h \in (0,1]$ it holds that:*

$$h^{n-1}\left|\Delta_h^n y\right| \leq K_0 \sup_{-(n-1)h \leq \tau \leq 0} |y(\tau)|, \forall y \in L_{loc}^\infty([-(n-1)h,0]; \Re) \qquad (2.4)$$

**Proof:** First notice that $P$ is a Vandermonde matrix and is invertible. Clearly, definition (2.2) guarantees that inequality (2.4) holds for appropriate $K_0 > 0$ (e.g., $K_0 = \sqrt{n}\left|P^{-1}\right|(n-1)!$). In order to show inequality (2.3), let $h \in (0,1]$, $u_i \in L_{loc}^\infty(\Re^+;\Re)$ ($i=1,...,n$), $x_0 \in \Re^n$ (arbitrary) and consider the solution $x(t)$ of (2.1) with initial condition $x(t_0) = x_0$ corresponding to inputs $u_i \in L_{loc}^\infty(\Re^+;\Re)$ ($i=1,...,n$). For all $t \geq t_0 + (n-1)h$, it holds that:

$$x(t) = \exp(A(n-1)h)x(t-(n-1)h) + \int_{t-(n-1)h}^{t} \exp(A(t-\tau))u(\tau)d\tau \qquad (2.5)$$

$$x_1(t-kh) = c'\exp(A(n-1-k)h)x(t-(n-1)h) + \int_{t-(n-1)h}^{t-kh} c'\exp(A(t-kh-\tau))u(\tau)d\tau, \; k=0,...,n-1 \qquad (2.6)$$

where $A = \begin{bmatrix} 0 & 1 & 0 & \cdots & 0 \\ 0 & 0 & 1 & \cdots & 0 \\ \vdots & \vdots & \vdots & & \vdots \\ 0 & 0 & 0 & \cdots & 1 \\ 0 & 0 & 0 & \cdots & 0 \end{bmatrix}$ and $c' = \begin{bmatrix} 1 & 0 & 0 & \cdots & 0 \end{bmatrix}$. Eliminating $x(t-(n-1)h)$ from the above equations, we obtain:

$$x(t) = Q_h^n P^{-1} \begin{bmatrix} x_1(t) \\ x_1(t-h) \\ x_1(t-2h) \\ \vdots \\ x_1(t-(n-1)h) \end{bmatrix} - Q_h^n P^{-1} B(t,u) + \int_{t-(n-1)h}^{t} \exp(A(t-\tau))u(\tau)d\tau, \; \forall t \geq t_0 + (n-1)h \qquad (2.7)$$

where

$$B'(t,u) = \left[ \int_{t-(n-1)h}^{t} c'\exp(A(t-\tau))u(\tau)d\tau \;\; \int_{t-(n-1)h}^{t-h} c'\exp(A(t-h-\tau))u(\tau)d\tau \;\; \cdots \;\; \int_{t-(n-1)h}^{t-(n-2)h} c'\exp(A(t-(n-2)h-\tau))u(\tau)d\tau \;\; 0 \right]$$

Since $c'\exp(Ah) = \begin{bmatrix} 1 & h & \frac{h^2}{2} & \cdots & \frac{h^{n-1}}{(n-1)!} \end{bmatrix}$, we have

$$\left| \int_{t-(n-1)h}^{t-kh} c'\exp(A(t-kh-\tau))u(\tau)d\tau \right| \leq \sum_{j=1}^{n} \frac{(n-1-k)^j}{j!} h^j \sup_{t-(n-1)h \leq \tau \leq t} |u_j(\tau)|$$

Consequently, $|B(t,u)| \leq (n-1)\sum_{j=1}^{n} \frac{(n-1)^j}{j!} h^j \sup_{t-(n-1)h \leq \tau \leq t} |u_j(\tau)|$ and

$$\left|Q_h^n P^{-1} B(t,u)\right| \leq (n-1)\sqrt{n}\left|P^{-1}\right|(n-1)!\sum_{j=1}^{n} \frac{(n-1)^j}{j!} h^{j+1-n} \sup_{t-(n-1)h \leq \tau \leq t} |u_j(\tau)|, \; \forall t \geq t_0 + (n-1)h \qquad (2.8)$$

Moreover, we have:

$$\left| \int_{t-(n-1)h}^{t} \exp(A(t-\tau))u(\tau)d\tau \right| \leq (n-1)h \max_{0 \leq s \leq (n-1)h} |\exp(As)| \sup_{t-(n-1)h \leq \tau \leq t} |u(\tau)|$$

$$\leq (n-1) \max_{0 \leq s \leq n-1} |\exp(As)| \sum_{j=1}^{n} h^{j+1-n} \sup_{t-(n-1)h \leq \tau \leq t} |u_j(\tau)|$$

The above inequality, in conjunction with (2.7), (2.8) and definition (2.2), guarantees that inequality (2.3) holds for appropriate constants $K_1,...,K_n > 0$ (e.g., $K_j = (n-1)\sqrt{n}|P^{-1}|(n-1)!\frac{(n-1)^j}{j!} + (n-1) \max_{0 \leq s \leq n-1} |\exp(As)|$, $j=1,...,n$). The proof is complete. ◁

The following lemma is concerned with the stabilization properties of system (1.1) by means of state ordinary linear feedback. Its proof is trivial and is omitted.

**Lemma 2.2** *Consider system (1.1), where $D \subset \Re^l$ is compact and $a : D \times \Re^n \to \Re$ is a locally Lipschitz mapping that satisfies (1.2). There exists a vector $k \in \Re^n$ and constants $M_0, M_1,..., M_n > 0$, $\mu > 0$ such that for every $v \in L^\infty_{loc}(\Re^+; \Re^n)$, $d \in L^\infty_{loc}(\Re^+; D)$, $(t_0, x_0) \in \Re^+ \times \Re^n$ the solution $x(t)$ of the closed-loop system (1.1) with $u = k'x$ corresponding to inputs $v \in L^\infty_{loc}(\Re^+; \Re^n)$, $d \in L^\infty_{loc}(\Re^+; D)$, with initial condition $x(t_0) = x_0 \in \Re^n$ satisfies:*

$$|x(t)| \leq M_0 \exp(-\mu(t-t_0))|x_0| + \sum_{i=1}^{n} M_i \sup_{t_0 \leq \tau \leq t} \{\exp(-\mu(t-\tau))|v_i(\tau)|\}, \forall t \geq t_0 \quad (2.9)$$

We are now in a position to state the main results of the present work. The following proposition shows that the state estimate provided by the backward difference operator of Lemma 2.1 can be used for the robust exponential stabilization of system (1.1) for sufficiently small values of the time step $h$. The result of the following proposition provides explicit formulae that allow the designer to select appropriate values for the time step $h$, in contrast with Proposition 2 in [11]. Its proof is provided at the Appendix.

**Proposition 2.3** *Consider system (1.1), where $D \subset \Re^l$ is compact and $a : D \times \Re^n \to \Re$ is a locally Lipschitz mapping that satisfies (1.2). Let $k \in \Re^n$ and $M_0, M_1,..., M_n > 0$, $\mu > 0$, be the vector and the constants for which (2.9) holds. Moreover, let $K_0, K_1,..., K_n > 0$, $\Delta^n_h : L^\infty_{loc}([-(n-1)h, 0]; \Re) \to \Re^n$, be the constants and the family of linear operators (parameterized by $h \in (0,1]$), for which (2.3), (2.4) hold. Then for every $h \in (0,1]$ that satisfies*

$$\beta h K_n |k| \exp(\mu(n-1)h) < 1; \quad \frac{hK_n M_n \beta^2 |k|^2 \exp(\mu(n-1)h)}{(1-\beta hK_n |k| \exp(\mu(n-1)h))^2} < 1 \quad (2.10)$$

*there exist constants $Q_0, Q_1,...,Q_n, Q_e > 0$ such that for every $v \in L^\infty_{loc}(\Re^+; \Re^n)$, $e \in L^\infty_{loc}(\Re; \Re)$ and $x_0 \in C^0([-(n-1)h, 0]; \Re^n)$, the solution $x(t)$ of the closed-loop system (1.1) with*

$$u(t) = k'\Delta^n_h T_{(n-1)h}(t)(x_1 + e) \quad (2.11)$$

*corresponding to inputs $v \in L^\infty_{loc}(\Re^+; \Re^n)$, $e \in L^\infty_{loc}(\Re; \Re)$, $d \in L^\infty_{loc}(\Re^+; D)$, with initial condition $x(\theta) = x_0(\theta); \theta \in [-(n-1)h, 0]$ satisfies (1.4) for all $t \geq 0$.*

**Remark 2.4:** Notice that formula (2.2) for the backward difference operator $\Delta^n_h : L^\infty_{loc}([-(n-1)h, 0]; \Re) \to \Re^n$ implies that the feedback law (2.11) can be equivalently expressed as

$$u(t) = k' Q_h^n P^{-1} \begin{bmatrix} x_1(t) + e(t) \\ x_1(t-h) + e(t-h) \\ \vdots \\ x_1(t-(n-1)h) + e(t-(n-1)h) \end{bmatrix}$$

where $Q_h^n := diag\left(1, \dfrac{1!}{(-h)}, \dfrac{2!}{(-h)^2}, \ldots, \dfrac{(n-1)!}{(-h)^{n-1}}\right) \in \Re^{n \times n}$ and $P = \begin{bmatrix} 1 & 0 & 0 & \ldots & 0 \\ 1 & 1 & 1^2 & \ldots & 1^{n-1} \\ 1 & 2 & 2^2 & \ldots & 2^{n-1} \\ \vdots & \vdots & \vdots & \ldots & \vdots \\ 1 & n-1 & (n-1)^2 & \ldots & (n-1)^{n-1} \end{bmatrix} \in \Re^{n \times n}$.

**Remark 2.5:** The proof of Proposition 2.3 (see Appendix) allows an estimation of the magnitude of the constants $Q_0, Q_1, \ldots, Q_n, Q_e > 0$ involved in (1.4):

$$Q_0 \leq M_0 + M_n \beta |k| \frac{cM_0 \beta |k| + nL}{1 - cM_n \beta^2 |k|^2}$$

$$Q_i \leq M_i + M_n \beta |k| \frac{c(K_i h^{i-n} + K_n \beta |k| M_i) + LK_n}{(1 - cM_n \beta^2 |k|^2) K_n}, \quad i = 1, \ldots, n$$

$$Q_e \leq M_n h^{1-n} K_0 \beta |k| \left[ \exp(\mu(n-1)h) + \beta |k| \frac{c(M_n \beta |k| \exp(\mu(n-1)h) + 1) + L \exp(\mu(n-1)h)}{1 - cM_n \beta^2 |k|^2} \right]$$

where $c := \dfrac{hK_n \exp(\mu(n-1)h)}{(1 - hK_n \beta |k| \exp(\mu(n-1)h))^2}$, $L := (1 + K_0 h^{1-n}) \exp(2(n + \beta |k| h^{1-n} K_0 + \mu)(n-1)h)$. However, numerical examples in the following section show that the above estimates are conservative for the linear case $a(d,x) \equiv 1$. Moreover, it should be noted that when $h \to 0$ then $Q_e \to +\infty$: this is the well-known phenomenon of sensitivity of high-gain feedback laws to measurement errors. Indeed, when $h \to 0$ then the linear feedback law (1.3) becomes a high-gain feedback (see the formula in the above remark).

Our next main result deals with the stabilization problem for system (1.5).

**Theorem 2.6** *Consider system (1.5), where $D \subset \Re^l$ is compact, the mappings $f: D \times \Re^k \times \Re^n \to \Re^k$, $g_i: D \times \Re^k \times \Re^n \to \Re$ ($i = 1, \ldots, n$) and $a: D \times \Re^k \times \Re^n \to \Re$ are locally Lipschitz with $f(d, 0, 0) = 0$ for all $d \in D$. Suppose the following:*

**(A1)** *There exist constants $\gamma \geq 0$, $c > 0$ and functions $a \in K_\infty$, $V \in C^0(\Re^k; \Re^+)$ which is positive definite and radially unbounded with $V(0) = 0$, such that for every $x \in L_{loc}^\infty(\Re^+; \Re^n)$, $d \in L_{loc}^\infty(\Re^+; D)$, $z_0 \in \Re^k$ the solution $z(t)$ of the system $\dot{z}(t) = f(d(t), z(t), x(t))$ corresponding to inputs $x \in L_{loc}^\infty(\Re^+; \Re^n)$, $d \in L_{loc}^\infty(\Re^+; D)$, with initial condition $z(0) = z_0 \in \Re^k$ satisfies:*

$$V(z(t)) \leq \exp(-ct) a(|z_0|) + \gamma \sup_{0 \leq \tau \leq t} \exp(-c(t-\tau)) |x(\tau)| \quad (2.12)$$

**(A2)** *There exists a constant $L \geq 0$, such that the following inequalities hold for all $(d, z, x) \in D \times \Re^k \times \Re^n$:*

$$|g_i(d, z, x)| \leq L |(x_1, \ldots, x_i)|, \quad i = 1, \ldots, n-1 \quad (2.13a)$$

$$|g_n(d, z, x)| \leq L(V(z) + |x|) \quad (2.13b)$$

**(A3)** *There exist constants $\alpha, \beta > 0$, such that the following inequalities hold for all $(d, z, x) \in D \times \mathfrak{R}^k \times \mathfrak{R}^n$:*

$$\alpha \leq a(d, z, x) \leq \beta \quad (2.13c)$$

Let $K_0, K_1, \ldots, K_n > 0$, $\Delta_h^n : L_{loc}^\infty([-(n-1)h, 0]; \mathfrak{R}) \to \mathfrak{R}^n$, be the constants and the family of linear operators (parameterized by $h \in (0,1]$), for which (2.3), (2.4) hold and let $k \in \mathfrak{R}^n$ and $M_0, M_1, \ldots, M_n > 0$, $\mu > 0$, be the vector and the constants for which (2.9) holds. Then for every $b \in (0,1]$ that satisfies

$$\frac{bK_n M_n \beta^2 |k|^2 \exp(\mu(n-1)b)}{(1 - bK_n \beta |k| \exp(\mu(n-1)b))^2} < 1; \quad bK_n \beta |k| \exp(\mu(n-1)b) < 1 \quad (2.14)$$

there exists $R(b) \geq 1$ and constants $M, K, Q > 0$ such that for every $r > R(b)$, $d \in L_{loc}^\infty(\mathfrak{R}^+; D)$, $v \in L_{loc}^\infty(\mathfrak{R}^+; \mathfrak{R}^n)$, $e \in L_{loc}^\infty(\mathfrak{R}; \mathfrak{R})$, $z_0 \in \mathfrak{R}^k$ and $x_0 \in C^0([-(n-1)h, 0]; \mathfrak{R}^n)$, where $h := b/r$, the solution $(z(t), x(t))$ of the closed-loop system (1.5) with

$$u(t) = r^n k' \, \text{diag}\!\left(1, \frac{1!}{(-b)}, \frac{2!}{(-b)^2}, \ldots, \frac{(n-1)!}{(-b)^{n-1}}\right) \begin{bmatrix} 1 & 0 & 0 & \ldots & 0 \\ 1 & 1 & 1^2 & \ldots & 1^{n-1} \\ 1 & 2 & 2^2 & \ldots & 2^{n-1} \\ \vdots & \vdots & \vdots & \ldots & \vdots \\ 1 & n-1 & (n-1)^2 & \ldots & (n-1)^{n-1} \end{bmatrix}^{-1} \begin{bmatrix} x_1(t) + e(t) \\ x_1(t-h) + e(t-h) \\ \vdots \\ x_1(t-(n-1)h) + e(t-(n-1)h) \end{bmatrix} \quad (2.15)$$

with initial condition $z(0) = z_0 \in \mathfrak{R}^k$, $x(\theta) = x_0(\theta); \theta \in [-(n-1)h, 0]$, corresponding to inputs $d \in L_{loc}^\infty(\mathfrak{R}^+; D)$, $v \in L_{loc}^\infty(\mathfrak{R}^+; \mathfrak{R}^n)$, $e \in L_{loc}^\infty(\mathfrak{R}; \mathfrak{R})$, satisfies for all $t \geq 0$:

$$|x(t)| + V(z(t)) \leq \exp(-\widetilde{\mu} r t) Q \left( p(r) \|x_0\|_{(n-1)h} + a(|z_0|) \right)$$
$$+ r^{-1} p(r) K \sup_{0 \leq \tau \leq t} \exp(-\widetilde{\mu} r(t-\tau)) |v(\tau)| + p(r) M \sup_{-(n-1)h \leq \tau \leq t} \exp(-\widetilde{\mu} r(t-\tau)) |e(\tau)| \quad (2.16)$$

where $\widetilde{\mu} := \min\!\left(\frac{c}{r}; \mu\right)$ and $p(r) = \frac{r^n}{r + 1 - R(b)}$.

**Remark 2.7:** Hypothesis (A1) is automatically satisfied if there exist constants $K, c, p > 0$ and a positive definite, radially unbounded and continuously differentiable function $W : \mathfrak{R}^k \to \mathfrak{R}^+$ such that for every $(d, z, x) \in D \times \mathfrak{R}^k \times \mathfrak{R}^n$ the differential inequality $\nabla W(z) f(d, z, x) \leq -2cpW(z) + K|x|^p$ holds. Particularly, in this case inequality (2.12) holds with $V(z) := (W(z))^{\frac{1}{p}}$, $\gamma := \left(\frac{K}{cp}\right)^{\frac{1}{p}}$ and $a \in K_\infty$ any function that satisfies $(W(z))^{\frac{1}{p}} \leq a(|z|)$ for all $z \in \mathfrak{R}^k$. Particularly, if the $z$-subsystem is linear $\dot{z} = Az + Bx$, where $B \in \mathfrak{R}^{k \times n}$ and $A \in \mathfrak{R}^{k \times k}$ is Hurwitz, then hypothesis (A1) holds for the function $W(z) = z'Pz$, where $P \in \mathfrak{R}^{k \times k}$ is a symmetric positive definite matrix which satisfies the property that the matrix $Q = -(A'P + PA)$ is positive definite.

**Proof of Theorem 2.6:** Notice that, by virtue of Remark 2.4, for every $d \in L_{loc}^\infty(\mathfrak{R}^+; D)$, $v \in L_{loc}^\infty(\mathfrak{R}^+; \mathfrak{R}^n)$, $e \in L_{loc}^\infty(\mathfrak{R}; \mathfrak{R})$, $z_0 \in \mathfrak{R}^k$ and $x_0 \in C^0([-(n-1)h, 0]; \mathfrak{R}^n)$, the solution $(z(t), x(t))$ of the closed-loop system (1.5) with (2.15) with initial condition $z(0) = z_0 \in \Re^k$, $x(\theta) = x_0(\theta); \theta \in [-(n-1)h, 0]$, corresponding to inputs $d \in L^\infty_{loc}(\Re^+; D)$, $v \in L^\infty_{loc}(\Re^+; \Re^n)$, $e \in L^\infty_{loc}(\Re; \Re)$, is related with the solution of

$$\frac{d\xi}{d\tau}(\tau) = r^{-1} f(\tilde{d}(\tau), \xi(\tau), diag(r, ..., r^n) w(\tau))$$

$$\frac{dw_i}{d\tau}(\tau) = r^{-(i+1)} g_i(\tilde{d}(\tau), \xi(\tau), diag(r, ..., r^n) w(\tau)) + w_{i+1}(\tau) + r^{-(i+1)} \tilde{v}_i(\tau) \quad, \quad i = 1, ..., n-1 \quad (2.17)$$

$$\frac{dw_n}{d\tau}(\tau) = r^{-(n+1)} g_n(\tilde{d}(\tau), \xi(\tau), diag(r, ..., r^n) w(\tau)) + a(\tilde{d}(\tau), \xi(\tau), diag(r, ..., r^n) w(\tau)) u(\tau) + r^{-(n+1)} \tilde{v}_n(\tau)$$

with

$$u(\tau) = k' \Delta_b^n T_{(n-1)b}(\tau)(w_1 + \tilde{e}) \quad (2.18)$$

initial condition $\xi(0) = z_0 \in \Re^k$, $w(\theta) = diag(r^{-1}, r^{-2}, ..., r^{-n}) x_0(\theta / r); \theta \in [-(n-1)b, 0]$, corresponding to inputs $\tilde{d} \in L^\infty_{loc}(\Re^+; D)$, $\tilde{v} \in L^\infty_{loc}(\Re^+; \Re^n)$, $\tilde{e} \in L^\infty_{loc}(\Re; \Re)$, where $\tilde{d}(\tau) := d(\tau / r)$, $\tilde{v}(\tau) := v(\tau / r)$, $\tilde{e}(\tau) := r^{-1} e(\tau / r)$, by the following formulae:

$$z(t) = \xi(rt) \; ; \; x(t) = diag(r, r^2, ..., r^n) w(rt), \text{ as long as the solutions exist} \quad (2.19)$$

By virtue of Proposition 2.3 and since (2.13c), (2.14) hold, there exist constants $Q_0, Q_1, ..., Q_n, Q_e > 0$ such that for every $\tilde{v} \in L^\infty_{loc}(\Re^+; \Re^n)$, $\tilde{e} \in L^\infty_{loc}(\Re; \Re)$ and $w_0 \in C^0([-(n-1)b, 0]; \Re^n)$, the solution $w(\tau)$ of the closed-loop system (2.17) with (2.18) corresponding to inputs $\tilde{v} \in L^\infty_{loc}(\Re^+; \Re^n)$, $\tilde{e} \in L^\infty_{loc}(\Re; \Re)$, with initial condition $w(\theta) = w_0(\theta); \theta \in [-(n-1)b, 0]$ satisfies the following inequality as long as the solution exists:

$$|w(\tau)| \le Q_0 \exp(-\mu \tau) \|w_0\|_{(n-1)b} + \sum_{i=1}^n Q_i r^{-(i+1)} \sup_{0 \le s \le \tau} \exp(-\mu(\tau - s)) |\tilde{v}_i(s)|$$

$$+ \sum_{i=1}^n Q_i r^{-(i+1)} \sup_{0 \le s \le \tau} \exp(-\mu(\tau - s)) |g_i(\tilde{d}(s), \xi(s), diag(r, ..., r^n) w(s))| + Q_e \sup_{-(n-1)b \le s \le \tau} \exp(-\mu(\tau - s)) |\tilde{e}(s)| \quad (2.20)$$

By virtue of (2.13a) we obtain for all $(\tilde{d}, \xi, w) \in D \times \Re^k \times \Re^n$:

$$r^{-(i+1)} |g_i(\tilde{d}, \xi, diag(r, ..., r^n) w)| \le \frac{L}{r} |w|, \; i = 1, ..., n-1 \quad (2.21)$$

Moreover, it follows from (2.13b) that the following inequality holds for all $(\tilde{d}, \xi, w) \in D \times \Re^k \times \Re^n$:

$$r^{-(n+1)} |g_n(\tilde{d}, \xi, diag(r, ..., r^n) w)| \le \frac{L}{r} |w| + \frac{L}{r^{n+1}} V(\xi), \; i = 1, ..., n-1 \quad (2.22)$$

It follows from (2.20), (2.21), (2.22) and the fact $\tilde{\mu} := \min\left(\frac{c}{r}; \mu\right) \le \mu$ that the following estimate holds as long as the solution of (2.17) with (2.18) exists:

$$\exp(\tilde{\mu} \tau) |w(\tau)| \le Q_0 \|w_0\|_{(n-1)b} + \sum_{i=1}^n Q_i r^{-(i+1)} \sup_{0 \le s \le \tau} \exp(\tilde{\mu} s) |\tilde{v}_i(s)| + \frac{L}{r^{n+1}} Q_n \sup_{0 \le s \le \tau} \exp(\tilde{\mu} s) V(\xi(s))$$

$$+ \frac{L}{r} \left( \sum_{i=1}^n Q_i \right) \sup_{0 \le s \le \tau} \exp(\tilde{\mu} s) |w(s)| + Q_e \sup_{-(n-1)b \le s \le \tau} \exp(\tilde{\mu} s) |\tilde{e}(s)| \quad (2.23)$$

On the other hand, hypothesis (A1) in conjunction with (2.19) guarantees that the following estimate holds as long as the solution of (2.17) with (2.18) exists:

$$V(\xi(\tau)) \leq \exp\left(-\frac{c}{r}\tau\right)a(|\xi(0)|) + \gamma r^n \sup_{0 \leq s \leq \tau} \exp\left(-\frac{c}{r}(\tau-s)\right)|w(s)| \qquad (2.24)$$

It follows from (2.24) and the fact $\tilde{\mu} := \min\left(\frac{c}{r}; \mu\right) \leq \frac{c}{r}$ that the following estimate holds as long as the solution of (2.17) with (2.18) exists:

$$\sup_{0 \leq s \leq \tau} \exp(\tilde{\mu} s)V(\xi(s)) \leq a(|\xi(0)|) + \gamma r^n \sup_{0 \leq s \leq \tau} \exp(\tilde{\mu} s)|w(s)| \qquad (2.25)$$

Combining (2.23) with (2.25), we conclude that the following estimate holds as long as the solution of (2.17) with (2.18) exists:

$$\sup_{0 \leq s \leq \tau} \exp(\tilde{\mu} s)|w(s)| \leq Q_0 \|w_0\|_{(n-1)b} + \sum_{i=1}^{n} Q_i r^{-(i+1)} \sup_{0 \leq s \leq \tau} \exp(\tilde{\mu} s)|\tilde{v}_i(s)| + \frac{L}{r^{n+1}} Q_n a(|\xi(0)|)$$
$$+ \frac{L}{r}\left(Q_n \gamma + \sum_{i=1}^{n} Q_i\right) \sup_{0 \leq s \leq \tau} \exp(\tilde{\mu} s)|w(s)| + Q_e \sup_{-(n-1)b \leq s \leq \tau} \exp(\tilde{\mu} s)|\tilde{e}(s)| \qquad (2.26)$$

Define:

$$R := R(b) = 1 + L\left(Q_n \gamma + \sum_{i=1}^{n} Q_i\right) \qquad (2.27)$$

Estimates (2.25), (2.26) and definition (2.27) in conjunction with the fact $r > R(b)$ give the following estimates holds as long as the solution of (2.17) with (2.18) exists:

$$\sup_{0 \leq s \leq \tau} \exp(\tilde{\mu} s)|w(s)| \leq \frac{r}{r+1-R} Q_0 \|w_0\|_{(n-1)b} + \sum_{i=1}^{n} Q_i \frac{r^{-i}}{r+1-R} \sup_{0 \leq s \leq \tau} \exp(\tilde{\mu} s)|\tilde{v}_i(s)| + \frac{L}{r^n(r+1-R)} Q_n a(|\xi(0)|)$$
$$+ Q_e \frac{r}{r+1-R} \sup_{-(n-1)b \leq s \leq \tau} \exp(\tilde{\mu} s)|\tilde{e}(s)| \qquad (2.28)$$

$$\sup_{0 \leq s \leq \tau} \exp(\tilde{\mu} s)V(\xi(s)) \leq \left(1 + \frac{L\gamma Q_n}{r+1-R}\right)a(|\xi(0)|) + \gamma \frac{r^{n+1}}{r+1-R} Q_0 \|w_0\|_{(n-1)b} + \sum_{i=1}^{n} \gamma Q_i \frac{r^{n-i}}{r+1-R} \sup_{0 \leq s \leq \tau} \exp(\tilde{\mu} s)|\tilde{v}_i(s)|$$
$$+ Q_e \gamma \frac{r^{n+1}}{r+1-R} \sup_{-(n-1)b \leq s \leq \tau} \exp(\tilde{\mu} s)|\tilde{e}(s)| \qquad (2.29)$$

Estimates (2.28) and (2.29) in conjunction with the fact that imply that the function $V \in C^0(\mathfrak{R}^k; \mathfrak{R}^+)$ is radially unbounded, imply that the phenomenon of finite escape time cannot happen. Thus the solution of (2.17) with (2.18) exists for all $\tau \geq 0$ and satisfies (2.28), (2.29) for all $\tau \geq 0$. By virtue of (2.19), it follows that the solution of the closed-loop system (1.5) with (2.15) exists for all $t \geq 0$. Exploiting (2.19) and definitions $\tilde{v}(\tau) := v(\tau/r)$, $\tilde{e}(\tau) := r^{-1}e(\tau/r)$ in conjunction with estimates (2.28), (2.29) and the facts $r > R(b) \geq 1$, $h := b/r$, gives for all $t \geq 0$:

$$\exp(\tilde{\mu} rt)|x(t)| + \exp(\tilde{\mu} rt)V(z(t)) \leq (1+\gamma)\frac{r^n}{r+1-R} Q_0 \|x_0\|_{(n-1)h} + \left(1 + \frac{(1+\gamma)LQ_n}{r+1-R}\right)a(|z_0|)$$
$$+ (1+\gamma)\sum_{i=1}^{n} Q_i \frac{r^{n-i}}{r+1-R} \sup_{0 \leq s \leq t} \exp(\tilde{\mu} rs)|v_i(s)| + Q_e(1+\gamma)\frac{r^n}{r+1-R} \sup_{-(n-1)h \leq s \leq t} \exp(\tilde{\mu} rs)|e(s)| \qquad (2.30)$$

Estimate (2.30) implies estimate (2.16) with $Q := (1+\gamma)Q_0 + 1 + (1+\gamma)LQ_n$, $K := (1+\gamma)\sum_{i=1}^{n} Q_i$ and $M := Q_e(1+\gamma)$.

Notice that the constants $M, K, Q > 0$ are all independent of $r$ (but depend on $b$). The proof is complete. ◁

## 3. Illustrating Examples

The following example illustrates the use of Proposition 2.3 for stabilization of linear systems.

**Example 3.1:** Consider the system (chain of three integrators)

$$\begin{aligned}\dot{x}_1 &= x_2 \\ \dot{x}_2 &= x_3 + v_2 \\ \dot{x}_3 &= u + v_3 \\ x &= (x_1, x_2, x_3)' \in \Re^3, u \in \Re, v_2 \in \Re, v_3 \in \Re\end{aligned} \quad (3.1)$$

The qualitative result in [11] guarantees that there exist $k_1, k_2, k_3$ and sufficiently small $h > 0$ such that the equilibrium point $0 \in C^0([-2h, 0]; \Re^3)$ is globally asymptotically stable for the closed-loop system (3.1) with $v_2 = v_3 \equiv 0$ and

$$u(t) = -k_1 x_1(t) - k_2 \frac{3x_1(t) - 4x_1(t-h) + x_1(t-2h)}{2h} - k_3 \frac{x_1(t) - 2x_1(t-h) + x_1(t-2h)}{h^2} \quad (3.2)$$

The use of Proposition 2.3 allows us to estimate the maximum allowable time step $h > 0$ that guarantees inequality (1.4) for the closed-loop system (3.1) with (3.2). Indeed, inequalities (2.3) and (2.4) hold for the operator

$$\Delta_h^3 y = \frac{1}{2h^2} \begin{bmatrix} 2h^2 y(0) \\ h(3y(0) - 4y(-h) + y(-2h)) \\ 2(y(0) - 2y(-h) + y(-2h)) \end{bmatrix} \quad (3.3)$$

with $K_0 \leq 4\sqrt{3}$ and $K_3 \leq \sqrt{136}$. Clearly, inequality (1.2) holds with $\alpha = \beta = 1$. Moreover, the vector $k = (-k_1, -k_2, -k_3)' = (-3, -5, -3)'$ guarantees that inequality (2.9) holds for the solution of the closed-loop system (1.1) with $n = 3$ and $u = k'x$. Particularly, inequality (2.9) holds with $\mu = 1/4$, $M_0 \leq \sqrt{190}$ and $M_3 \leq 2\sqrt{5}$. The estimation of the constants $\mu$, $M_0$ and $M_3$ is performed by making use of the quadratic Lyapunov function $V(x) = \frac{1}{2}x_1^2 + \frac{1}{2}(x_2 + x_1)^2 + \frac{1}{2}(x_3 + 2x_2 + 2x_1)^2$. It follows that inequalities (2.10) hold for $h \leq 4.1 \cdot 10^{-4}$. However, (as noted above in Remark 2.5) it should be emphasized that this is a conservative estimate of the maximum allowable time step. Numerical simulations have shown that the maximum allowable time step is approximately $h = 0.21$. Of course, as $h \to 0.21$, the rate of convergence becomes slower ($\mu \to 0$ in (1.4)). In Figure 1 it is shown the evolution of the states of the closed-loop system (3.1) with (3.2), $h = 0.1$ and $v_2 = v_3 \equiv 0$ (initial condition $x_2(0) = x_3(0) = 1$, $x_1(\theta) = 0$ for $\theta \in [-0.2, -0.1]$ and $x_1(\theta) = 10\theta + 1$ for $\theta \in [-0.1, 0]$). Clearly, the states converge to zero exponentially, despite the fact that the estimates provided by the backward difference operator (3.3) are not good approximations of the state vector during the interval $[0, 0.2]$. The initial transient period, where the estimates provided by the backward difference operator (3.3) are not good approximations of the state vector, deteriorates the performance (compare with Figure 2, which shows the evolution of the states of the closed-loop system (3.1) with $u = k'x$, $v_2 = v_3 \equiv 0$ and initial condition $x_1(0) = x_2(0) = x_3(0) = 1$; as expected state feedback guarantees better performance compared to output feedback with $x_1$ as output).

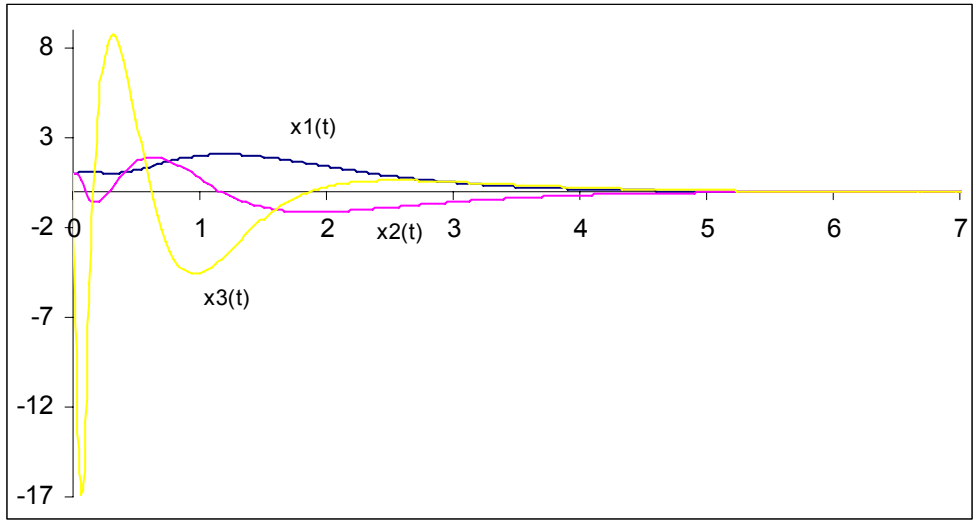

**Figure 1:** The evolution of the states of the closed-loop system (3.1) with (3.2), $h = 0.1$ and $v_2 = v_3 \equiv 0$

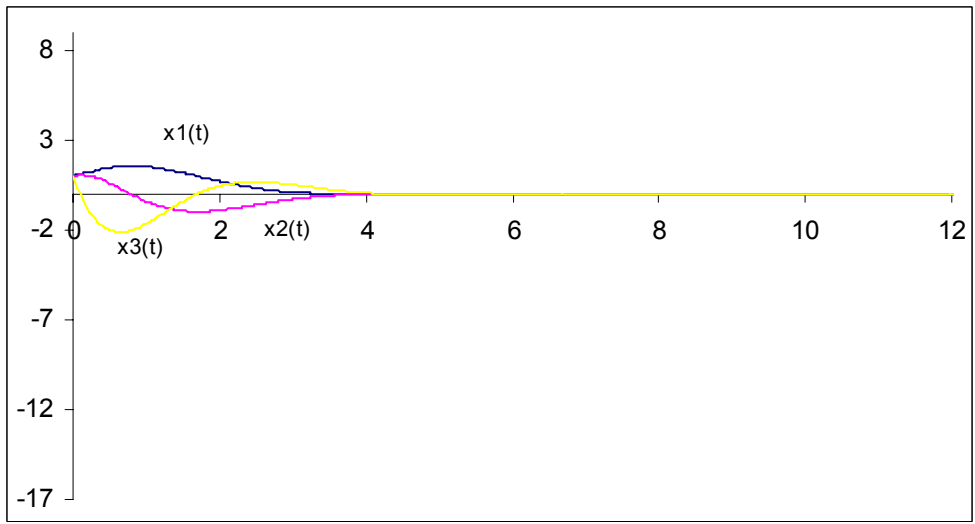

**Figure 2:** The evolution of the states of the closed-loop system (3.1) with $u = k'x$, $v_2 = v_3 \equiv 0$

In Figure 3 it is shown the evolution of the states of the closed-loop system (3.1) with (3.2), $h = 0.1$ and $v_2(t) = \cos(t)$, $v_3(t) = 1.5\sin(t)$ (initial condition $x_2(0) = x_3(0) = 1$, $x_1(\theta) = 0$ for $\theta \in [-0.2, -0.1]$ and $x_1(\theta) = 10\theta + 1$ for $\theta \in [-0.1, 0]$). Clearly, the states converge to a periodic solution exponentially and estimate (1.4) holds.

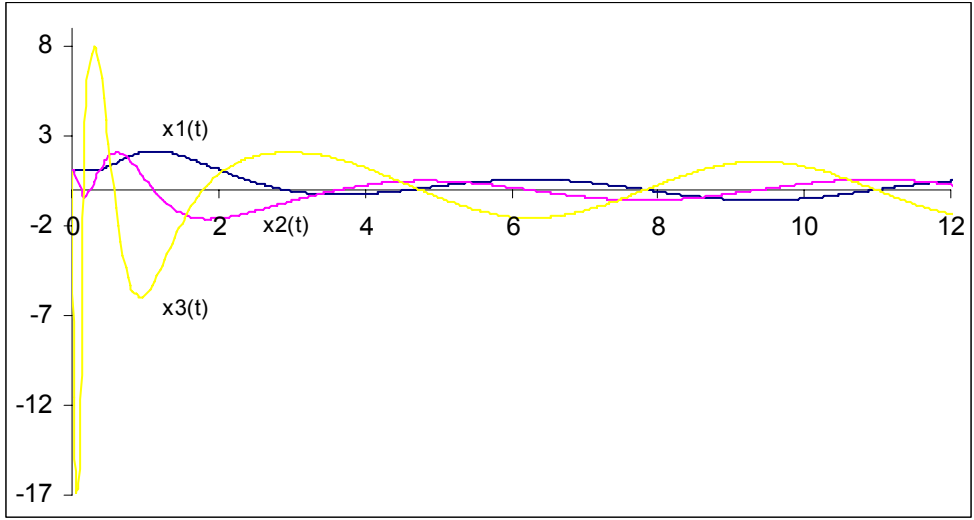

**Figure 3:** The evolution of the states of the closed-loop system (3.1) with (3.2), $h = 0.1$ and $v_2(t) = \cos(t)$, $v_3(t) = 1.5\sin(t)$

**Example 3.2:** Consider the following nonlinear system:

$$\dot{z} = -z - z^3 + d_1 x_2$$
$$\dot{x}_1 = x_2 \; ; \; \dot{x}_2 = x_3 \; ; \; \dot{x}_3 = d_2 z^2 + u \qquad (3.4)$$
$$z \in \Re, x = (x_1, x_2, x_3)' \in \Re^3, u \in \Re, d = (d_1, d_2) \in [-1,1] \times [-1,1]$$

By virtue of Remark 2.7 and using the function $W(z) = \frac{1}{4} z^4$, it follows that hypothesis (A1) of Theorem 2.6 holds with $V(z) := \frac{1}{2} z^2$, $\gamma := \frac{1}{2}$ and $a(s) := \frac{1}{4} s^4$. Hypotheses (A2), (A3) of Theorem 2.6 hold as well with $L := 2$, $\alpha = \beta = 1$. Consequently, Theorem 2.6 guarantees that there exist $k_1, k_2, k_3$ and sufficiently small $h > 0$ such that the equilibrium point $0 \in C^0([-2h,0];\Re^4)$ is robustly globally asymptotically stable for the closed-loop system (3.4) with (3.2). In Figures 4 and 5, it is shown the evolution of the states of the closed-loop system (3.4) with (3.2), $h = 0.1$, $k = (-k_1, -k_2, -k_3)' = (-3, -5, -3)'$, $d_1(t) = \text{sgn}(x_2(t))$, $d_2(t) \equiv 1$ (initial condition $z(0) = 2$, $x_2(0) = x_3(0) = 1$, $x_1(\theta) = 0$ for $\theta \in [-0.2, -0.1]$ and $x_1(\theta) = 10\theta + 1$ for $\theta \in [-0.1, 0]$).

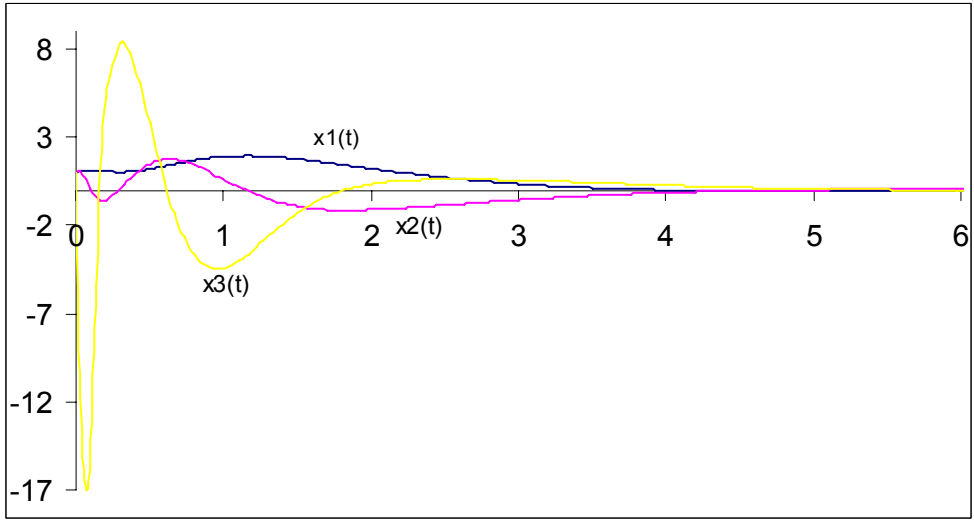

**Figure 4:** The evolution of the state $x(t)$ of the closed-loop system (3.4) with (3.2), $h = 0.1$

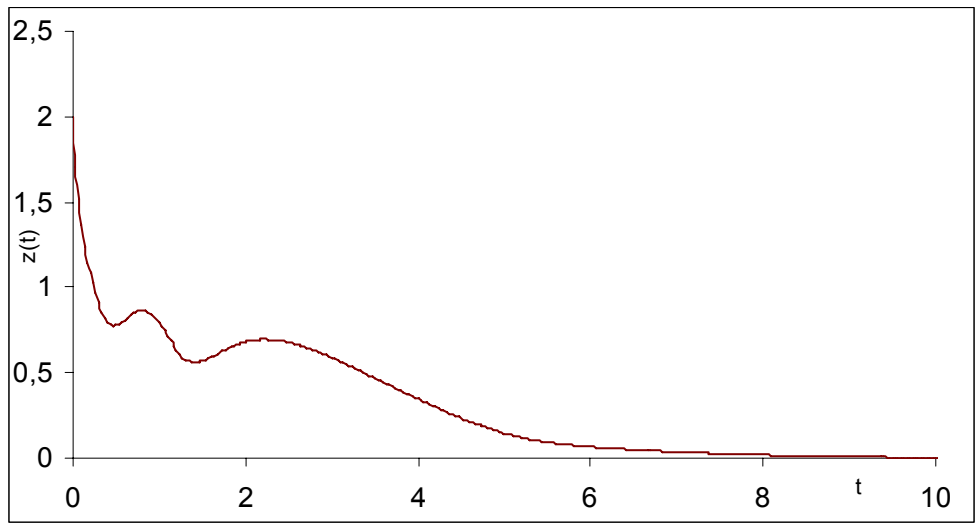

**Figure 5:** The evolution of the state $z(t)$ of the closed-loop system (3.4) with (3.2), $h = 0.1$

## 4. Concluding Remarks

In this work, sufficient conditions for robust global asymptotic stabilization of nonlinear uncertain systems by means of discrete-delay static output feedback are presented. The efficiency of the proposed control strategy is demonstrated with illustrating examples. It is clear that static output discrete-delay feedback is an alternative to dynamic output (observer-based) feedback and further studies need to be performed, in order to show the advantages and disadvantages of each type of feedback.

# Appendix

**Proof of Proposition 2.3:** For notational convenience we set $r = (n-1)h$. Following the proof of Theorem 1.1 (page 168) in [2] and using (2.4), we notice that the solution $x(t)$ of the closed-loop system (1.1) with (2.11) corresponding to arbitrary inputs $v \in L_{loc}^\infty(\Re^+;\Re^n)$, $e \in L_{loc}^\infty(\Re;\Re)$, with initial condition $x(\theta) = x_0(\theta); \theta \in [-r,0]$ exists for all $t \geq 0$ and satisfies the estimate:

$$\sup_{t-r \leq \tau \leq t} |x(\tau)| \leq \exp(M\,t)\left(n\|x_0\|_r + \sum_{i=1}^n \sup_{0 \leq \tau \leq t}|v_i(\tau)| + \beta|k|h^{1-n}K_0 \sup_{-r \leq \tau \leq t}|e(\tau)|\right), \quad \forall t \geq 0 \quad (A1)$$

where (by virtue of (1.2), (2.4)) $M = n + \beta|k|h^{1-n}K_0$. Furthermore, by virtue of (1.2), (2.3), (2.4), we have for all $t \geq r$:

$$\left|x(t) - \Delta_h^n T_r(t)x_1\right|$$

$$\leq \sum_{j=1}^{n-1} K_j h^{j+1-n} \sup_{t-r \leq \tau \leq t}|v_j(\tau)| + K_n h \sup_{t-r \leq \tau \leq t}\left|v_n(\tau) + a(d(\tau), x(\tau))\left(k'x(\tau) - k'\left(x(\tau) - \Delta_h^n T_r(\tau)x_1\right) + k'\Delta_h^n T_r(\tau)e\right)\right| \quad (A2)$$

$$\leq \sum_{j=1}^n K_j h^{j+1-n} \sup_{t-r \leq \tau \leq t}|v_j(\tau)| + \beta h K_n|k|\left(\sup_{t-r \leq \tau \leq t}|x(\tau)| + \sup_{t-r \leq \tau \leq t}\left|x(\tau) - \Delta_h^n T_r(\tau)x_1\right| + h^{1-n}K_0 \sup_{t-r \leq \tau \leq t}|e(\tau)|\right)$$

It follows from (A2) that the following estimate holds for all $t \geq 2r$:

$$\sup_{t-r \leq \tau \leq t}\exp(\mu\tau)\left|x(\tau) - \Delta_h^n T_r(\tau)x_1\right| \leq \sum_{j=1}^n K_j h^{j+1-n}\exp(\mu r)\sup_{t-2r \leq \tau \leq t}\exp(\mu\tau)|v_j(\tau)|$$

$$+ \beta h K_n|k|\exp(\mu r)\left(\sup_{t-2r \leq \tau \leq t}\exp(\mu\tau)|x(\tau)| + \sup_{t-2r \leq \tau \leq t}\exp(\mu\tau)\left|x(\tau) - \Delta_h^n T_r(\tau)x_1\right| + h^{1-n}K_0 \sup_{t-2r \leq \tau \leq t}\exp(\mu\tau)|e(\tau)|\right)$$

By distinguishing the cases $\sup_{t-2r \leq \tau \leq t}\exp(\mu\tau)\left|x(\tau) - \Delta_h^n T_r(\tau)x_1\right| = \sup_{t-2r \leq \tau \leq t-r}\exp(\mu\tau)\left|x(\tau) - \Delta_h^n T_r(\tau)x_1\right|$ and $\sup_{t-2r \leq \tau \leq t}\exp(\mu\tau)\left|x(\tau) - \Delta_h^n T_r(\tau)x_1\right| = \sup_{t-r \leq \tau \leq t}\exp(\mu\tau)\left|x(\tau) - \Delta_h^n T_r(\tau)x_1\right|$, it follows from (2.10) and the above inequality that the following estimate holds for all $t \geq 2r$:

$$\sup_{t-r \leq \tau \leq t}\exp(\mu\tau)\left|x(\tau) - \Delta_h^n T_r(\tau)x_1\right| \leq \sum_{j=1}^n \frac{K_j h^{j+1-n}\exp(\mu r)}{1 - \beta h K_n|k|\exp(\mu r)}\sup_{t-2r \leq \tau \leq t}\exp(\mu\tau)|v_j(\tau)|$$

$$+ \frac{\beta h K_n|k|\exp(\mu r)}{1 - \beta h K_n|k|\exp(\mu r)}\left(\sup_{t-2r \leq \tau \leq t}\exp(\mu\tau)|x(\tau)| + h^{1-n}K_0 \sup_{t-2r \leq \tau \leq t}\exp(\mu\tau)|e(\tau)|\right) \quad (A3)$$

$$+ \beta h K_n|k|\exp(\mu r)\sup_{t-2r \leq \tau \leq t-r}\exp(\mu\tau)\left|x(\tau) - \Delta_h^n T_r(\tau)x_1\right|$$

Using estimate (A3) and induction, it may be shown that the following inequality holds for all $m \in Z^+$, $\xi \geq 2r$:

$$\sup_{\xi+(m-1)r \leq \tau \leq \xi+mr}\exp(\mu\tau)\left|x(\tau) - \Delta_h^n T_r(\tau)x_1\right| \leq \sum_{j=1}^n \frac{K_j h^{j+1-n}\exp(\mu r)}{1-\beta h K_n|k|\exp(\mu r)}\sum_{l=0}^m (\beta h K_n|k|\exp(\mu r))^{m-l} \sup_{\xi+(l-2)r \leq \tau \leq \xi+lr}\exp(\mu\tau)|v_j(\tau)|$$

$$+ \frac{\beta h K_n|k|\exp(\mu r)}{1-\beta h K_n|k|\exp(\mu r)}\sum_{l=0}^m (\beta h K_n|k|\exp(\mu r))^{m-l}\left(\sup_{\xi+(l-2)r \leq \tau \leq \xi+lr}\exp(\mu\tau)|x(\tau)| + h^{1-n}K_0 \sup_{\xi+(l-2)r \leq \tau \leq \xi+lr}\exp(\mu\tau)|e(\tau)|\right)$$

$$+ (\beta h K_n|k|\exp(\mu r))^{m+1}\sup_{\xi-2r \leq \tau \leq \xi-r}\exp(\mu\tau)\left|x(\tau) - \Delta_h^n T_r(\tau)x_1\right|$$

(A4)

Let arbitrary $T \geq 2r$. Using inequality (A4) with $m = \left[\dfrac{T-2r}{r}\right]$, $\xi = t - r\left[\dfrac{T-2r}{r}\right]$ (where $\left[\dfrac{T-2r}{r}\right]$ denotes the integer part of the non-negative real number $\dfrac{T-2r}{r}$), in conjunction with the fact that $\sum\limits_{l=0}^{m}(\beta h K_n |k| \exp(\mu r))^{m-l} \leq \dfrac{1}{1-\beta h K_n |k| \exp(\mu r)}$ for all $m \in Z^+$, we obtain:

$$\sup_{T-r \leq \tau \leq T} \exp(\mu \tau)\left|x(\tau) - \Delta_h^n T_r(\tau) x_1\right| \leq \sum_{j=1}^{n} \frac{K_j h^{j+1-n} \exp(\mu r)}{(1 - \beta h K_n |k| \exp(\mu r))^2} \sup_{0 \leq \tau \leq T} \exp(\mu \tau)|v_j(\tau)|$$
$$+ \frac{\beta h K_n |k| \exp(\mu r)}{(1 - \beta h K_n |k| \exp(\mu r))^2}\left(\sup_{0 \leq \tau \leq T} \exp(\mu \tau)|x(\tau)| + h^{1-n} K_0 \sup_{0 \leq \tau \leq T} \exp(\mu \tau)|e(\tau)|\right)$$
$$+ \exp(-\sigma(T-2r)) \sup_{0 \leq \tau \leq 2r} \exp(\mu \tau)\left|x(\tau) - \Delta_h^n T_r(\tau) x_1\right|$$

where $\sigma := -\dfrac{1}{r}\log(\beta h K_n |k|) - \mu > 0$. Clearly,

$$\sup_{r \leq \tau \leq t} \exp(\mu \tau)\left|x(\tau) - \Delta_h^n T_r(\tau) x_1\right| \leq \sup_{2r \leq T \leq t}\left(\sup_{T-r \leq \tau \leq T} \exp(\mu \tau)\left|x(\tau) - \Delta_h^n T_r(\tau) x_1\right|\right), \text{ for all } t \geq 2r$$

Hence, the two above inequalities give for all $t \geq 2r$:

$$\sup_{r \leq \tau \leq t} \exp(\mu \tau)\left|x(\tau) - \Delta_h^n T_r(\tau) x_1\right| \leq \sum_{j=1}^{n} \frac{K_j h^{j+1-n} \exp(\mu r)}{(1 - \beta h K_n |k| \exp(\mu r))^2} \sup_{0 \leq \tau \leq t} \exp(\mu \tau)|v_j(\tau)|$$
$$+ \frac{\beta h K_n |k| \exp(\mu r)}{(1 - \beta h K_n |k| \exp(\mu r))^2}\left(\sup_{0 \leq \tau \leq t} \exp(\mu \tau)|x(\tau)| + h^{1-n} K_0 \sup_{0 \leq \tau \leq t} \exp(\mu \tau)|e(\tau)|\right) \qquad (A5)$$
$$+ \sup_{0 \leq \tau \leq 2r} \exp(\mu \tau)\left|x(\tau) - \Delta_h^n T_r(\tau) x_1\right|$$

On the other hand, it follows from (1.2), (2.9) and (2.4), that the solution $x(t)$ of the closed-loop system (1.1) with (2.11) corresponding to inputs $v \in L_{loc}^{\infty}(\Re^+; \Re^n)$, $e \in L_{loc}^{\infty}(\Re; \Re)$, with initial condition $\{x(\theta); \theta \in [-r,0]\} = x_0 \in C^0([-r,0]; \Re^n)$, satisfies the following estimate for all $t \geq 0$:

$$\exp(\mu t)|x(t)| \leq M_0|x(0)| + \sum_{i=1}^{n} M_i \sup_{0 \leq \tau \leq t}\{\exp(\mu \tau)|v_i(\tau)|\}$$
$$+ \beta M_n |k| \sup_{0 \leq \tau \leq t}\left\{\exp(\mu \tau)\left|x(\tau) - \Delta_h^n T_r(\tau) x_1\right|\right\} + \beta M_n |k| h^{1-n} K_0 \exp(\mu r) \sup_{-r \leq \tau \leq t} \exp(\mu \tau)|e(\tau)| \qquad (A6)$$

Combining estimates (A5) and (A6), we obtain for all $t \geq 2r$:

$$\sup_{r \leq \tau \leq t} \exp(\mu \tau)\left|x(\tau) - \Delta_h^n T_r(\tau) x_1\right| \leq \sum_{j=1}^{n} \frac{h \exp(\mu r)(K_j h^{j-n} + \beta K_n |k| M_j)}{(1 - \beta h K_n |k| \exp(\mu r))^2} \sup_{0 \leq \tau \leq t} \exp(\mu \tau)|v_j(\tau)|$$
$$+ \frac{\beta h K_n |k| \exp(\mu r)}{(1 - \beta h K_n |k| \exp(\mu r))^2} M_0|x(0)| + \frac{\beta h K_n |k| \exp(\mu r)}{(1 - \beta h K_n |k| \exp(\mu r))^2} \beta M_n |k| \sup_{0 \leq \tau \leq t}\left\{\exp(\mu \tau)\left|x(\tau) - \Delta_h^n T_r(\tau) x_1\right|\right\} \qquad (A7)$$
$$+ \sup_{0 \leq \tau \leq 2r} \exp(\mu \tau)\left|x(\tau) - \Delta_h^n T_r(\tau) x_1\right| + \frac{\beta h^{2-n} K_n |k| \exp(\mu r)}{(1 - \beta h K_n |k| \exp(\mu r))^2} K_0(\beta M_n |k| \exp(\mu r) + 1) \sup_{-r \leq \tau \leq t} \exp(\mu \tau)|e(\tau)|$$

By using (2.10) and distinguishing the cases $\sup_{0 \leq \tau \leq t} \{\exp(\mu\tau)|x(\tau) - \Delta_h^n T_r(\tau)x_1|\} = \sup_{0 \leq \tau \leq r} \{\exp(\mu\tau)|x(\tau) - \Delta_h^n T_r(\tau)x_1|\}$ and $\sup_{0 \leq \tau \leq t} \{\exp(\mu\tau)|x(\tau) - \Delta_h^n T_r(\tau)x_1|\} = \sup_{r \leq \tau \leq t} \{\exp(\mu\tau)|x(\tau) - \Delta_h^n T_r(\tau)x_1|\}$, inequality (A7) implies for all $t \geq 2r$ :

$$\sup_{r \leq \tau \leq t} \exp(\mu\tau)|x(\tau) - \Delta_h^n T_r(\tau)x_1| \leq \sum_{j=1}^{n} \frac{c(K_j h^{j-n} + \beta K_n |k| M_j)}{(1 - cM_n \beta^2 |k|^2) K_n} \sup_{0 \leq \tau \leq t} \exp(\mu\tau)|v_j(\tau)|$$
$$+ \frac{cM_0 \beta |k|}{1 - cM_n \beta^2 |k|^2}|x(0)| + h^{1-n} \frac{cK_0 \beta |k|(M_n \beta |k| \exp(\mu r) + 1)}{1 - cM_n \beta^2 |k|^2} \sup_{-r \leq \tau \leq t} \exp(\mu\tau)|e(\tau)| \quad (A8)$$
$$+ \frac{1}{1 - cM_n \beta^2 |k|^2} \sup_{0 \leq \tau \leq 2r} \exp(\mu\tau)|x(\tau) - \Delta_h^n T_r(\tau)x_1|$$

where $c := \dfrac{hK_n \exp(\mu r)}{(1 - \beta h K_n |k| \exp(\mu r))^2}$. Notice that by virtue of (2.4) and (A1), we have for all $t \leq 2r$ :

$$\sup_{0 \leq \tau \leq t} \exp(\mu\tau)|x(\tau) - \Delta_h^n T_r(\tau)x_1|$$
$$\leq (1 + K_0 h^{1-n}) \exp(2(n + \beta |k| h^{1-n} K_0 + \mu) r) \left( n \|x_0\|_r + \sum_{i=1}^{n} \sup_{0 \leq \tau \leq t} |v_i(\tau)| + \beta |k| h^{1-n} K_0 \sup_{-r \leq \tau \leq t} |e(\tau)| \right) \quad (A9)$$

Finally, combining estimates (A8) and (A9) we obtain for all $t \geq 0$ :

$$\sup_{0 \leq \tau \leq t} \exp(\mu\tau)|x(\tau) - \Delta_h^n T_r(\tau)x_1| \leq \sum_{j=1}^{n} \frac{c(K_j h^{j-n} + \beta K_n |k| M_j) + LK_n}{(1 - cM_n \beta^2 |k|^2) K_n} \sup_{0 \leq \tau \leq t} \exp(\mu\tau)|v_j(\tau)|$$
$$+ \frac{cM_0 \beta |k| + nL}{1 - cM_n \beta^2 |k|^2} \|x_0\|_r + h^{1-n} \beta K_0 |k| \frac{c(M_n \beta |k| \exp(\mu r) + 1) + L \exp(\mu r)}{1 - cM_n \beta^2 |k|^2} \sup_{-r \leq \tau \leq t} \exp(\mu\tau)|e(\tau)| \quad (A10)$$

where $L := (1 + K_0 h^{1-n}) \exp(2(n + \beta |k| h^{1-n} K_0 + \mu) r)$. Estimate (A10) in conjunction with estimate (A6) gives for all $t \geq 0$ :

$$\exp(\mu t)|x(t)| \leq \left( M_0 + M_n \beta |k| \frac{cM_0 \beta |k| + nL}{1 - cM_n \beta^2 |k|^2} \right) \|x_0\|_r$$
$$+ \sum_{i=1}^{n} \left[ M_i + M_n \beta |k| \frac{c(K_i h^{i-n} + K_n \beta |k| M_i) + LK_n}{(1 - cM_n \beta^2 |k|^2) K_n} \right] \sup_{0 \leq \tau \leq t} \exp(\mu\tau)|v_i(\tau)|$$
$$+ M_n h^{1-n} K_0 \beta |k| \left[ \exp(\mu r) + \beta |k| \frac{c(M_n \beta |k| \exp(\mu r) + 1) + L \exp(\mu r)}{1 - cM_n \beta^2 |k|^2} \right] \sup_{-r \leq \tau \leq t} \exp(\mu\tau)|e(\tau)|$$

The above estimate implies (1.4). The proof is complete. ◁